\input amstex

\input epsf
\mag=\magstep1
\documentstyle{amsppt}
\pagewidth{135mm}
\pageheight{200mm}
\vcorrection{-5mm}
\nopagenumbers

\document

\topmatter
\title An application of non-positively curved cubings of alternating links \endtitle
\author Makoto Sakuma and Yoshiyuki Yokota \endauthor
\thanks
\endthanks
\thanks
The first author was partially supported by JSPS KAKENHI Grant Number 15H03620.
The second author was partially supported by JSPS KAKENHI Grant Number 15K04878.
\endthanks
\abstract
By using non-positively curved cubings of prime alternating link exteriors,
we prove that certain ideal triangulations of their complements, 
derived from reduced alternating diagrams, are non-degenerate,
in the sense that none of the edges is homotopic relative its endpoints to a peripheral arc.
This guarantees that 
the hyperbolicity equations for those triangulations for hyperbolic alternating links
have solutions corresponding to the complete hyperbolic structures.
Since the ideal triangulations considered in this paper are often used in the study of the volume conjecture, 
this result has a potential application to the volume conjecture.
\endabstract
\endtopmatter

\head 1. Introduction \endhead

Let $L$ be a hyperbolic link in $S^3$ and $M:=S^3\setminus L$ the complement of $L$.
In the approach to the volume conjecture\,[{\bf 7},\,{\bf 8}],
initiated by D. Thurston\,[{\bf 9}] and the second author\,[{\bf 11},\,{\bf 12}],
a certain ideal triangulation, $\Cal S$, of $M$ derived from a diagram of $K$ plays a crucial role.
Tetrahedra in $\Cal S$ correspond to $q$-factorials in the Kashaev invariant,
and the hyperbolicity (gluing) equations for $\Cal S$ and the complex volume of $M$ are related to
the potential function which appears in an integral expression of the Kashaev invariant.

However, in general, there is no guarantee that the hyperbolicity equations for $\Cal S$ have a geometric solution, 
i.e., a solution which corresponds to the complete hyperbolic structure.
In fact, this fails if and only if some edge of $\Cal S$ is {\it inessential}
in the sense that it is homotopic to a curve on the peripheral torus.
To be precise, for a link $L$ in $S^3$, the arc in its exterior $E(L):=S^3\setminus N(L)$
obtained from the ideal edge of $\Cal S$ is homotopic, relative to its endpoints, to an arc in $\partial E(L)$.
(Here, $N(L)$ is an open regular neighborhood of $L$.)
In [{\bf 6}], such an arc is called {\it peripheral}.
An inessentail ideal edge has no geodesic representative in the hyperbolic manifold $M$,
and conversely, if all ideal edges of $\Cal S$ are {\it essential} (i.e., not inessential),
then the edges have unique geodesic representatives,
which gives a geometric solution to the hyperbolicity equations
(though some of the tetrahedra may be flat or negatively oriented).

The purpose of this paper is to prove the following theorem,
by using non-positively curved cubings of alternating link exteriors.

\proclaim{Theorem 1.1}
Let $L$ be a hyperbolic link in $S^3$ which has a reduced alternating diagram $D$
and $\Cal S$ the ideal triangulation of the complement of $L$ associated to $D$.
Then, the edges of $\Cal S$ are essential,
and so the hyperbolicity equations obtained from $\Cal S$ have a geometric solution.
\endproclaim

Non-positively curved cubings of alternating link exteriors were first found
by the pioneering work by Aitchison, Lumsden and Rubinstein\,[{\bf 3}]:
they proved that if an alternating link $L$ admits a \lq\lq nicely balanced'' alternating diagram,
then its exterior admits a non-positively curved cubing. 
The existence of non-positively curved cubing for the exterior of every prime alternating link
was first noted in a literature by Adams\,[{\bf 1}], where he attributes it to Agol.
In fact, Agol gave a beautiful application of non-positively curved cubings of the exteriors of $2$-bridge links,
in his talk\,[{\bf 2}] in 2002.

The cubings of link exteriors themselves have been essentially known to the experts from early time. 
Ideal triangulations derived from the related cubings are used in the wonderful computer program {\it SnapPea},
as explained by Weeks\,[{\bf 13}].
The cubing $\Cal C$ of a link exterior $E(L)$ is intimately related to the classical Dehn complex $\Cal D$ of $L$:
there is a deformation retraction $r:E(L)\to \Cal D$,
and $\Cal C$ is identified with the mapping cylinder of the restriction of $r$ to $\partial E(L)$\,(see Section 2).
It is a classical result due to Weinbaum\,[{\bf 14}]
that the Dehn presentation of the augmented link group $\pi_1(E(L))*\Bbb Z$ obtained from a diagram $D$ of $L$
satisfies the small cancellation condition if and only if $D$ is a reduced alternating diagram.
Moreover, it is known that the Dehn complex $\Cal D$ obtained from a link diagram $D$ is non-positively curved
if and only if $D$ is a prime alternating diagram\,(see [{\bf 4},\,Proposition I\!I.5.43]).
This is translated to the corresponding cubing as follows:
the cubed complex $\Cal C$ obtained from a link diagram $D$ is non-positively curved
if and only if $D$ is a prime alternating diagram\,(see Proposition 3.3 and Remark 3.4).

This paper is organized as follows.
In Section 2, we review the cubings and ideal triangulations of the exteriors and complements of prime alternating links.
In Section 3, we review some basic facts on cubed complexes, and apply to alternating link exteriors.
We give a proof of Theorem 1.1 in Section 4 and state its consequence in Section 5. 

\medskip
{\bf Acknowledgements.}
The authors realized that the results of this paper can be derived 
from the existence of non-positively curved cubings of alternating hyperbolic link exteriors in 2015,
and it was announced by the second author at a conference at Waseda University
in honour of the 20th anniversary of the Volume Conjecture.
During the conference, we learned from Stavros Garoufalidis that
essentially the same results had been already obtained in 2002 
by the joint work of S. Garoufalidis, I. Moffatt and D. Thurston in [{\bf 6}],
which was completed in 2007, but has not been published.
Their method is based on the small cancellation property of the Dehn presentation.
Though their proof is algebraic while ours is geometric,
both proofs are based on the non-positively curved property of alternating link exteriors.
We thank S. Garoufalidis for kindly sharing their preprint with us and encouraging us to publish our results.
We also thank the organizer, Jun Murakami, of the conference at Waseda University.

\head 2. Cubical decomposition of link exteriors \endhead

Let $L$ be a prime link in $S^3$ which is represented by a connected diagram $D$,
and $E(L):=S^3\setminus N(L)$ the {\it exterior} of $L$, where $N(L)$ is an open regular neighborhood of $L$.
D. Thurston\,[{\bf 9}] described a method for decomposing $E(L)$
into partially truncated octahedra placed between crossings.
He also described a method for constructing an ideal tetrahedral decomposition
of the link complement from the octahedral decomposition.
These decompositions are essentially equivalent to those described by Weeks\,[{\bf 13}],
and their details are described by the second author\,[{\bf 12}].
The octahedral decomposition induces a cubical decomposition of $E(L)$,
which is non-positively curved if and only if the diagram $D$ is reduced and alternating,
as observed in Remark 3.4.

In this section, we recall the cubical decomposition and the ideal triangulation obtained from it, following [{\bf 12}].
In what follows, we assume that $D$ is a reduced alternating diagram.
We may pick two points $P_+$ and $P_-$ in $S^3$, identify $S^3\setminus\{P_+,P_-\}$ with $S^2\times\Bbb R$,
and assume the following hold.
The diagram $D$ is regarded as a 4-valent graph in $S^2\times\{0\}$,
$L\subset D\times[-1,1]\subset S^2\times [-1,1]$,
and $L$ intersects $S^2\times\{0\}$ transversely in $2c$ points, $P_1,P_2,\dots,P_{2c}$,
where $c$ is the crossing number of $D$.
Let $D^*$ be a graph embedded in $S^2\times\{0\}$ dual to $D$, such that
$D\cap D^*=\{P_1,P_2,\dots,P_{2c}\}$.
The vertices of $D$ and $D^*$ are denoted by $X_1,X_2,\dots,X_c$ and $R_0,R_1,\dots,R_{c+1}$ respectively.
The closures of the $4c$ connected components of $S^2\setminus(D\cup D^*)$ are denoted by $Q_1,Q_2,\dots,Q_{4c}$ and
$$
\align
\nu&:\{1,2,\dots,4c\}\to\{1,2,\dots,c\},\\
\mu&:\{1,2,\dots,4c\}\to\{0,1,\dots,c+1\},\\
\alpha,\beta&:\{1,2,\dots,4c\}\to\{0,1,\dots,2c-1\},\\
\sigma&:\{1,2,\dots,4c\}\to\{-1,1\}
\endalign
$$
are defined by Figure 1, where $1\le g\le 4c$.
$$
\gather
\epsfxsize=200pt
\epsfbox{cube01}\\
\text{Figure 1}
\endgather
$$
The connected components of
$$
L_+=L\cap(S^2\times[0,1]),\quad L_-=L\cap(S^2\times[-1,0])
$$
are called overpasses and underpasses of $L$ respectively,
and we assume that each overpass/underpass intersects $S^2\times \{\pm 1\}$
precisely at the point above/below the corresponding crossing (vertex) of $D$. Namely, we assume
$$
L\cap(S^2\times\{\pm1\})=\{X_n\,|\,1\le n\le c\}\times\{\pm1\}.
$$
Observe that $S^2\times\Bbb R$ is decomposed into $4c$ quadratic prisms
$$
Q_1\times\Bbb R,\,Q_2\times\Bbb R,\,\dots,\,Q_{4c}\times\Bbb R
$$
each of which looks as in Figure 2,
where the bold arcs represent $I_g:=L_+\cap(Q_g\times\Bbb R)$ and $J_g:=L_-\cap(Q_g\times\Bbb R)$.
$$
\gather
\epsfxsize=160pt
\epsfbox{cube02}\\
\text{Figure 2}
\endgather
$$
We consider the arcs
$$
\gather
A(X_n):=X_n\times(1,\infty),\ B(X_n):=X_n\times(-\infty,-1),\ C(X_n):=X_n\times(-1,1),\\
A(P_l):=P_l\times(0,\infty),\ B(P_l):=P_l\times(-\infty,0),\ C(R_m):=R_m\times\Bbb R,
\endgather
$$
where $n\in\{1,2,\dots,c\}$, $m\in\{0,1,\dots,c+1\}$ and $l\in\{0,1,\dots,2c-1\}$.
Let $A(I_g)$ and $B(J_g)$ be the 2-cells in $\partial(Q_g\times\Bbb R)$ bounded by
$$
\{I_g,A(X_{\nu(g)}),A(P_{\alpha(g)})\},\quad\{J_g,B(X_{\nu(g)}),B(P_{\beta(g)})\}
$$
respectively, the shaded ones in Figure 2, where $1\le g\le 4c$.
Then, as each of the overpasses and underpasses is contractible,
we can collapse each connected component of $\cup_{g=1}^{4c}A(I_g)$ and $\cup_{g=1}^{4c}B(J_g)$ to a vertical edge.
Then each $Q_g\times\Bbb R$ becomes an ideal tetrahedron, $S_g$, in $(S^2\times\Bbb R)\setminus L$ with 6 ideal edges
$$
A(X_{\nu(g)})=A(P_{\alpha(g)}),\ A(P_{\beta(g)}),\ 
B(X_{\nu(g)})=B(P_{\beta(g)}),\ B(P_{\alpha(g)}),\ C(X_{\nu(g)}),\ C(R_{\mu(g)})
$$
as shown in Figure 3. 
The family of ideal tetrahedra, $\{S_g\}_{1\le g\le 4c}$,  
gives an ideal triangulation, $\Cal T$, of $(S^2\times\Bbb R)\setminus L$.
$$
\gather
\epsfxsize=200pt
\epsfbox{cube03}\\
\text{Figure 3}
\endgather
$$
Note that $\{Q_g\times\Bbb R:\nu(g)=n\}$ and $\{S_g:\nu(g)=n\}$ intersect
$\partial N(P_\pm),\partial N(L)$ as shown in Figures 4a and 4b,
where $N(P_\pm)$ denote open regular neighborhoods of $P_\pm$ respectively.
$$
\gather
\epsfxsize=56pt
\epsfbox{cube04a}\hskip15mm
\epsfxsize=56pt
\epsfbox{cube04b}\hskip 15mm
\epsfxsize=56pt
\epsfbox{cube04c}\\
\text{Figure 4a}\hskip19mm\text{Figure 4b}\hskip20mm\text{Figure 4c}
\endgather
$$
For each $n\in\{1,2,\dots,c\}$, the union $O_n:=\cup_{\nu(g)=n}S_g$ of the $4$ ideal tetrahedra
sharing the \lq\lq crossing arc'' $C(X_n)$ is regarded as a quotient of an ideal octahedron.
To see this, note that $S_g$ is identified with the join $C(X_{\nu(g)})*C(R_{\mu(g)})$.
Thus $O_n$ is identified with the join of $C(X_n)$ and $Z_n:=\cup_{\nu(g)=n}C(R_{\mu(g)})$, which is a cycle of length $4$.
This implies that $O_n$ is regarded as a quotient of an ideal octahedron.
In fact, we can obtain an ideal octahedron by cutting $O_n$
along $A(X_n)$ and $B(X_n)$\,(see [{\bf 9}, Figure in p.17] and [{\bf 11}, Figure 2]).
In the following, we do not distinguish between $O_n$ and an ideal octahedron.

By adding the two vertices $P_{\pm}$ to $O_n$
(i.e., by \lq\lq replacing" each of the four ideal vertices of the octahedron contained in the cycle $Z_n$ with a real vertex),
and taking the intersection with $E(L)$
(i.e., truncating the pair of ideal vertices corresponding to the two ends of the crossing arc $C(X_n)$),
we obtain a partially truncated octahedron, $\check O_n$, and the set $\{\check O_n\}_{1\le n\le c}$ determines
a partially truncated octahedral decomposition of $E(L)$, where $P_{\pm}$ are the only vertices in the interior of $E(L)$.

Now note that the join of the midpoint $X_n$ of the crossing arc $C(X_n)$
with the cycle $Z_n$ determines an ideal square in $O_n$, which divides $O_n$ into two ideal pyramids.
The ideal square in $O_n$ descends to a square in $\check O_n$,
and the pair of pyramids descends to a pair of cubes in $\check O_n$
which intersect $\partial N(P_\pm)$ and $\partial N(L)$ as shown in Figure 4c.
The set of these pairs of cubes determines a cubic decomposition, $\Cal C$, of $E(L)$.
The cubic decomposition has precisely two vertices $P_{\pm}$ in the interior of $E(L)$ and the vertices
$$
a_n:=A(X_n)\cap\partial N(L),\ b_n:=B(X_n)\cap\partial N(L)
$$
in $\partial E(L)$, where $1\le n\le c$. 

We note that the set of the squares determines the Dehn complex, $\Cal D$, of $L$.
For each cube in $\Cal C$, there is a unique square in $\Cal D$ which is a face of the cube.
Thus there is a natural deformation retraction $r:\Cal C=E(L)\to\Cal D$
and $\Cal C$ is identified with the mapping cylinder of the restriction of $r$ to $\partial E(L)$.

At the end of this section, we recall the ideal triangulation, $\Cal S$, of $E(L)$
which are used in the study of the volume conjecture by [{\bf 9},\,{\bf 11},\,{\bf 12}].
Roughly speaking,
$\Cal S$ is obtained from the ideal triangulation $\Cal T$ of $S^2\times\Bbb R=S^3\setminus (L\cup\{P_+,P_-\})$
by engulfing the extra ideal points $P_{\pm}$ into $L$, where we suppose $c>3$\,(see Assumption 3 in [{\bf 12}]).
To this end, pick a point, say $P_0$, from $L\cap S^2\times\{0\}$,
and collapse the ideal edges $A(P_0)$ and $B(P_0)$ into an ideal vertex.
This forces some of the other simplexes of $\Cal T$ to degenerate.
In order to give more precise description,
arrange $Q_1,Q_2,Q_{4c-1},Q_{4c},X_1,X_c,R_0,R_{c+1},P_1,P_{2c-1}$ as shown in Figure 5,
namely the following hold for the functions introduced in this section.
$$
\gather
\nu(1)=\nu(2)=\nu(\alpha^{-1}(1))=1,\ \nu(4c)=\nu(4c-1)=\nu(\beta^{-1}(2c-1))=c,\\
\mu(1)=\mu(4c-1)=0,\ \mu(2)=\mu(4c)=c+1,\\
\alpha(1)=\alpha(2)=\beta(4c-1)=\beta(4c)=0\\
\epsfxsize=150pt
\epsfbox{cube05}\\
\text{Figure 5}
\endgather
$$
Then, since $A(P_0)=A(X_1)=A(P_1)$, $B(P_0)=B(X_c)=B(P_{2c-1})$,
and since the collapsing of $A(P_0)$ and $B(P_0)$ cause collapsing of $C(R_0)$ and $C(R_{c+1})$ into ideal vertices,
the following degeneration of ideal tetrahedra occur.
The pairs of ideal tetrahedra
$\{S_1,S_2\}$ and $\{S_{4c-1},S_{4c}\}$ collapse into the ideal edges $B(X_1)$ and $A(X_c)$ respectively,
and, if $g$ belongs to
$$
(\alpha^{-1}(\{1,2c-1\})\cup\beta^{-1}(\{1,2c-1\})\cup\mu^{-1}(0,c+1))\setminus\{1,2,4c-1,4c\},
$$
then $S_g$ collapses into an ideal triangle.
The other ideal tetrahedra continue to be ideal tetrahedra, namely, if $g$ belongs to
$$
\Gamma=\{1,2,\dots,4c\}\setminus(\alpha^{-1}(\{1,2c-1\})\cup\beta^{-1}(\{1,2c-1\})\cup\mu^{-1}(0,c+1)),
$$
then $S_g$ remains to be an ideal tetrahedron in $\Cal S$. See [{\bf 12}] for details.
(Though [{\bf 12}] treats only hyperbolic knot diagrams which satisfy certain assumptions,
the same arguments are available for prime link diagrams which satisfy the same assumptions,
such as reduced alternating diagrams with $\ge4$ crossings.)

In particular, the edges of $\Cal S$ are obtained as the images of the following paths in $\Cal C$.
$$
\alignat 2
\alpha_n:&=A(X_1)\cup P_+\cup A(X_n) & (1<n\le c)&\\
\beta_n:&=B(X_c)\cup P_-\cup B(X_n)  & (1\le n<c)&\\
\gamma_n:&=C(X_n) & (1\le n\le c)&\\
\delta_m:&=A(X_1)\cup P_+\cup C(R_m)\cup P_-\cup B(X_c) & \qquad(m\not\in\mu(\nu^{-1}(\{1,c\})))&
\endalignat
$$

\head 3. Non-positively curved cubing \endhead

Let $L$, $D$ and $\Cal C$ be as in the previous section.
Then, by identifying each cube in $\Cal C$ with a unit cube in $\Bbb E^3$,
$\Cal C$ is regarded as a {\it cubed complex}\,(see [{\bf 4}, p.115]).
In this section, we explain the well-known fact that $\Cal C$ is non-positively curved.

We first review some basic facts about a cubed complex $W$.
A path $\gamma:[0,l]\to W$ is called a {\it piecewise geodesic} if there exist
$$
0=t_0<\cdots<t_{i-1}<t_i<\cdots<t_n=l
$$
such that each $\gamma|_{[t_{i-1},t_i]}$ is an isometric embedding into some cube.
Then there exists a shortest piecewise geodesic between any two points in $X$,
called a {\it geodesic}, and $W$ becomes a complete geodesic metric space.

For $x\in W$, the set of unit tangent vectors at $x$ is called the {\it geometric link} of $x$ in $W$
and denoted by $\text{Lk}(x,W)$.
We can regard $\text{Lk}(x,W)$ as a piecewise spherical complex,
and so it admits the structure of a complete geodesic metric space. 
Note that the length of each edge of $\text{Lk}(x,W)$ is equal to $\pi/2$.
The following lemma is well-known\,(see [{\bf 4}, Remark I.5.7]).

\proclaim{Lemma 3.1}
A piecewise geodesic $\gamma$ in $W$ is a local geodesic if, for each point $x$ on $\gamma$,
the distance between the incoming and the outgoing unit tangent vectors to $\gamma$ at $x$ are
at least $\pi$ in $\text{\rm Lk}(x,W)$.
\endproclaim

In general, a metric space is said to be {\it non-positively curved} if each point in it has a neighborhood
where any geodesic triangle is thinner than a comparison triangle in $\Bbb E^2$,
that is, the distance between any points on a geodesic triangle is less than or equal to
the distance between the corresponding points on a comparison triangle\,(cf. [{\bf 4}, Definition I\!I.1.2]).
The following criterion is well-known\,(see [{\bf 4}, Theorem I\!I.5.20]).

\proclaim{Proposition 3.2(Gromov's link condition)} A cubed complex $W$ is non-positively curved
if $\text{\rm Lk}(x,W)$ is a simplicial complex which is flag,
i.e., any finite subset of vertices, that is pairwisely joined by edges, spans a simplex.
\endproclaim

The following proposition was first noted in a literature by Adams\,[{\bf 1}], where he attributes it to Agol.

\proclaim{Proposition 3.3}
Let $L$ be a prime alternating link in $S^3$ represented by a reduced alternating diagram $D$,
and let $\Cal C$ be the cubed complex with underlying space $E(L)$ constructed from the diagram $D$.
Then $\Cal C$ and its double, $\hat\Cal C$, across $\partial E(L)$ are non-positively curved.
\endproclaim

\demo{proof}
Observe that $\Cal C$ induces a cubing of $\partial E(L)$, such that each vertex has degree $4$\,(see Figure 4c).
Thus, for any vertex of $\Cal C$ contained in $\partial E(L)$,
its link in $\Cal C$ is identified with a unit hemi-sphere 
consisting of $4$ spherical triangles which are regular and right-angled.
Hence these vertices satisfy Gromov's link condition.

For the inner vertex $P_{+}$ of $\Cal C$,
we can observe that the link $\text{Lk}(P_+,\Cal C)$ is obtained from the cell decomposition of $S^2$
determined by the graph $D^*$, by subdiving each region of $D^*$ as follows. 
Each region of $D^*$ contains a unique vertex, say $X_n$, of $D$.
Subdivide the region by taking the join of $X_n$
and the edge cycle of $D^*$ forming the boundary of the region\,(see Figure 6).
In fact, the vertex $X_n$ of $\text{Lk}(P_+,\Cal C)$ comes from the edge $A(X_n)$ of $\Cal C$,
whereas the vertex $R_m (\in D^*)$ of $\text{Lk}(P_+,\Cal C)$ comes from the edge $C(R_m)$ of $\Cal C$\,(see Figure 4c).

The checker board coloring of the regions of $D$ induces 
a black and white coloring of the vertices of $D^*$.
Thus the set of the vertices of $\text{Lk}(P_+,\Cal C)$
is divided into the following three subsets: 
the set of the white vertices of $D^*$, 
the set of the black vertices of $D^*$,
and the set of the vertices of $D$.
Moreover, any edge of $\text{Lk}(P_+,\Cal C)$
joins vertices which belong to different groups.
Now pick a triple of vertices of $\text{Lk}(P_+,\Cal C)$ such that
each sub-pair spans an edge in $\text{Lk}(P_+,\Cal C)$.
Then the triple contains a unique vertex from each subset.
In particular, it contains a pair of vertices which form the boundary of an edge, say $e$, of $D^*$. 
The only vertices of $\text{Lk}(P_+,\Cal C)$ which span an edge with each of the boundary vertices of $e$
are the two vertices of $D$ dual to the two regions of $D^*$ containing $e$ in the boundary.
Thus the triple of the vertices span an $2$-simplex of $\text{Lk}(P_+,\Cal C)$.
Hence $\text{Lk}(P_+,\Cal C)$ satisfies Gromov's condition.
The same argument works for $\text{Lk}(P_-,\Cal C)$.
Thus we have proved that $\Cal C$ is non-positively curved.
$$
\gather
\epsfxsize=220pt
\epsfbox{cube06}\\
\text{Figure 6}
\endgather
$$

We can easily check that the double $\hat\Cal C$ is also non-positively curved,
because the link in $\hat\Cal C$ of a vertex contained in $\partial E(L)$
is the double of the link in $\Cal C$ of the corresponding vertex 
and so it is identified with the unit $2$-sphere consisting of
8 spherical triangles which are regular and right-angled. \qed
\enddemo

\proclaim{Remark 3.4}
{\rm The cubing $\Cal C$ can be constructed from any connected link diagram $D$,
and we can see as in the above proof that $\Cal C$ is non-positively curved
if and only if $D$ is a reduced alternating diagram.}
\endproclaim

\head 4. Proof of Theorem 1.1\endhead

We show that the arcs $\alpha_n, \beta_n, \gamma_n, \delta_m$, 
introduced at the end of Section 2, are local geodesics in the cubed complex $\Cal C$.
Then, since these arcs are orthogonal to the boundary $\partial\Cal C$,
this implies that their doubles in the double $\hat\Cal C$ of $\Cal C$ are closed local geodesics.
Since $\hat\Cal C$ is non-positively curved by Proposition 3.3,
this in turn implies that these loops are not null-homotopic in $\hat\Cal C$\,(see [{\bf 4}, Theorem I\!I.4.23]).
Hence the arcs $\alpha_n, \beta_n, \gamma_n, \delta_m$ are essential\,(non-peripheral) in $\Cal C=E(L)$.
Since these arcs form the edge set of $\Cal S$, we obtain the desired result.

We first show that $\alpha_n=A(X_1)\cup P_+\cup A(X_n)$ with $1<n\le c$ is a local geodesic in $\Cal C$.
Note that $\alpha_n$ is a piecewise geodesic consisting of two geodesic arcs $A(X_1)$ and $A(X_n)$,
where the two arcs \lq\lq intersect" $\text{Lk}(P_+,\Cal C)$ at the vertices $X_1$ and $X_n$.
Since these two vertices belong to the same subset introduced in the proof of Proposition 3.3,
no edge of $\text{Lk}(P_+,\Cal C)$ joins them.
Since $\text{Lk}(P_+,\Cal C)$ is a spherical complex  consisting of right-angled regular triangles,
this implies that the distance between the vertices $X_1$ and $X_n$ in $\text{Lk}(P_+,\Cal C)$ is $\ge \pi$.
Hence $\alpha_n$ is a local geodesic by Lemma 3.1. The same argument works for $\beta_n$.

Next, show that $\delta_m=A(X_1)\cup P_+\cup C(R_m)\cup P_-\cup B(X_c)$
with $m\not\in\mu(\nu^{-1}(\{1,c\}))$ is a local geodesic in $\Cal C$.
We have only to show that $\delta_m$ satisfies the condition in Lemma 3.1 at the vertices $P_+$ and $P_-$.
To this end, note that the intersection of $\delta_m$ with $\text{Lk}(P_+,\Cal C)$ are the vertices $X_1$ and $R_m$.
Since $m\not\in\mu(\nu^{-1}(\{1\}))$, no edge of $\text{Lk}(P_+,\Cal C)$ joins them,
and hence the distance between $X_1$ and $R_m$ in $\text{Lk}(P_+,\Cal C)$ is $\ge \pi$.
Thus $\delta_m$ satisfies the condition in Lemma 3.1 at $P_+$.
The same argument works for $P_-$. Hence $\delta_m$ is a local geodesic.

The remaining arcs $\gamma_n$ are obviously local geodesics,
and so we have proved that the arcs $\alpha_n, \beta_n, \gamma_n, \delta_m$ are local geodesics.
This completes the proof of Theorem 1.1.
The fact that the crossing edges $\gamma_n$ are essential are already known by [{\bf 1}] and [{\bf 5}],
as is noted in [{\bf 10}].

\head 5. Application \endhead

Let $L$ be a hyperbolic link in $S^3$ with a reduced alternating diagram $D$ and $M$ its complement.
Under the same assumption in Section 2, the potential function associated to $D$ is defined by 
$$
V(z;D)=\sum_{g\in\Gamma}\sigma(g)\text{Li}_2(z(\beta(g))/z(\alpha(g))),
$$
where $z$ is a map $\{2,3,\dots,2c-2\}\to\Bbb C$ such that $z(l)=1$ if $l\in\alpha(\mu^{-1}(\{0,c+1\}))$.
For simplicity, we put
$$
\Lambda=\{2,3,\dots,2c-2\}\setminus\alpha(\mu^{-1}(\{0,c+1\})).
$$
Then, by Theorem 1.1 and [{\bf 12},\,Theorems 2.5 and 2.6], we have the following corollary.
(The results in [{\bf 12}] holds not only for hyperbolic knots but also )

\proclaim{Corollary 5.1} The hyperbolicity equations for $\Cal S$ are given by
$$
\exp\left\{z(l)\frac{\partial V(z;D)}{\partial z(l)}\right\}=1,\quad l\in\Lambda
$$
which must have a solution $\zeta$ corresponding to the complete hyperbolic structure of $M$. Furthermore,
$$
\hat V(\zeta,D)=V(\zeta;D)-\sum_{l\in\Lambda}\log\zeta(l)\left[z(l)\frac{\partial V(z;D)}{\partial z(l)}\right]_{z=\zeta}
$$
gives the complex volume of $M$ modulo $\pi^2$.
\endproclaim

At the end of this paper, we note that Ian Agol\,[{\bf 2}] had announced
the following beautiful application of non-positively curved cubings of hyperbolic $2$-bridge link exteriors.
{\it Each hyperbolic $2$-bridge link group admit precisely two parabolic gerenating pairs,
namely the upper and lower meridian pairs.}
In fact, he has shown that any pair of parabolic transformations, 
which are not equivalent to the upper and lower meridian pairs,
generates a free subgroup of the $2$-bridge link group.
The proof is based on the non-positivity of the cubing
and the fact that the union of the \lq\lq vertical middle planes" of the cubes in $\Cal C$ 
give rise to the checker board surfaces associated with the reduced alternating diagram.
(The union of the horizontal middle planes gives the boundary of a regular neighborhood of the Dehn complex.)
The same argument implies that any pair of parabolic transformations
of the link group of a hyperbolic alternating link with bridge index $\ge 3$
generates a free subgroup of the link group.

\vskip 5mm\noindent
{\eightrm Department of Mathematics, Guraduate School of Science, Hiroshima University,}
\par\noindent
{\eightrm Higashi-Hiroshima, 739-8526, Japan}
\par\noindent
{\eightrm E-mail address: sakuma\@hiroshima-u.ac.jp}

\vskip 5mm\noindent
{\eightrm Department of Mathematics and Information Sciences, Tokyo Metropolitan University,}
\par\noindent
{\eightrm Hachioji, Tokyo, 192-0397, Japan}
\par\noindent
{\eightrm E-mail address: jojo\@tmu.ac.jp}

\enddocument